\documentclass[12pt]{amsart}
\usepackage[centertags]{amsmath}
\usepackage{amsfonts}
\usepackage{amssymb}
\usepackage{amsthm}
\usepackage{newlfont}
\usepackage{bm}
\newlength{\defbaselineskip}
\newcommand{\setlinespacing}[1]%
{\setlength{\baselineskip}{#1 \defbaselineskip}}

\def\DD{{\mathcal D}}

\def\KK{{\mathcal K}}
\def\LL{{\mathcal L}}
\def\MM{{\mathcal M}}

\def\TT{{\mathcal T}}

\def\bbC{\mathbb{C}}

\def\bbD{\mathbb{D}}
\def\bbT{\mathbb{T}}

\def\1{\mathbf{1}}

\def\bbC{\mathbb{C}}

\def\bbD{\mathbb{D}}
\def\bbT{\mathbb{T}}

\theoremstyle{plain}
\newtheorem{thm}{Theorem}[section]

\newtheorem{lem}[thm]{Lemma}
\newtheorem{prop}[thm]{Proposition}

\theoremstyle{definition}
\newtheorem{defn}{Definition}[section]
\def\dss{\displaystyle}
\theoremstyle{remark}

\numberwithin{equation}{section}

\begin{document}

\def\dss{\displaystyle}

\title{Product of Matrix Valued Truncated Toeplitz Operators }

\author{Muhammad Ahsan Khan}
\address{Department of Mathematics University of Sialkot, Pakistan}
\email{muhammad.ahsan@uskt.edu.pk}
\begin{abstract}
Let $A_\Phi$ be a matrix valued truncated Toeplitz operator-the compression of multipication operator to vector valued model space $H^2(E)\ominus \Theta H^2(E)$, where $\Theta$ is a matrix valued non constant inner function. Under supplementary assumptions we find necessary and sufficient condition that the product $A_\Phi A_\Psi$ is itself a matrix valued truncated Toeplitz operator. 	
\end{abstract}
\maketitle
\section{Introduction}
	 
Toeplitz operators are the compressions of multiplication operator to the usual Hardy Hilbert space $H^2$. In \cite{BH}, Brown and Halmos describe the algebraic properties of Toeplitz operators. Among other things, they found necessary and sufficient conditions for the product of two Toeplitz operators to itself be a Toeplitz operator, namely that either the first operator’s symbol is antiholomorphic or the second operator’s symbol is holomorphic. In either case, the symbol of the product is the product of the symbols.\\
In the last decade, a large amount of research has concentrated on a generalization
of Toeplitz matrices, namely truncated Toeplitz operators. These are the compressions of multiplication operator to subspaces of the Hardy space which are invariant under the backward shift operator. They have been formally introduced in \cite{Sa}; see \cite{GR} for a more recent survey.  Sarason ~\cite{Sa} found equivalents to several of Brown and Halmos’s results for truncated Toeplitz operators on the model spaces $H^2\ominus \theta H^2$, where $\theta$ is some non-constant inner function. The model spaces are the backward-shift invariant subspaces of $H^2$
(that they are backward shift invariant follows easily from the fact that $\theta H^2$ is clearly shift invariant). We refer the reader to ~\cite{BBK,BCFMT,GR} (see also \cite{GMR}) for the general theory of these operators. It is well know that the product of two truncated Toeplitz operators is not a truncated Toeplitz operator. In particular, in \cite{NS} Sedlock has investigated when a product of truncated Toeplitz operators is itself a truncated Toeplitz operator.\\ Most recently, the basics of corresponding matrix valued truncated Toeplitz operators (MTTOs), which are compressions to $\KK_\Theta$ of multiplications with matrix valued functions on $H^2(E)$ has been developed in \cite{KT}. In view of the result of \cite{NS}, there arises a basic question related to the product of matrix valued truncated Toeplitz operators that is  when the product of two MTTOs is still an MTTO? But there is no such simple result in the general case, and we need some supplementary assumptions to obtain the main result Theorem 4.6. The purpose of the present paper is to adapt the approach in \cite{CT} to the case of MTTOs on an arbitrary model space. \\
The plan of the paper is following:
By means of Section~2, we want to make sure that the reader has become acquainted to model spaces and their operators and other useful facts from this area, needed when we are going to start the main work in section~4. In Section ~3 we will introduce truncated Toeplitz operators (TTOs) and matrix valued truncated  Toeplitz operators (MTTOs). The last section contains a particular case of MTTOs namely the Block Toeplitz matrices.
\section{Model Spaces and Operators}  
Let $\bbC$ denote the complex plane, $\bbD$ the unit disc in $\bbC$ and $\bbT$ one dimensional torus in $\mathbb{C}$. In the sequel $E$ will denote a fixed Hilbert space of dimension $d$. We designate the algebra of bounded linear operators on $E$ by $\mathcal{L}(E)$ and by $\LL(E,K)$ the space of all bounded linear operators from Hilbert space $E$ to a Hilbert space $K$. The space $L^2(E)$ is defined as usual, by 
	 \[
	 L^2(E):=\left\{f:\bbT\longrightarrow E: f(e^{it})=\sum_{n=-\infty}^{\infty} a_ne^{int}: a_n\in E,\sum_{n=-\infty}^{\infty}\|a_n\|^2<\infty\right\}
	 \] 
	 endowed with the inner product
	 \[\langle f , g \rangle_{L^2 (E)} =\int_{\bbT}\langle f(e^{it}), g(e^{it})\rangle _E dm
	 \]
	 where $dm=\frac{dt}{2\pi}$ is the normalized \emph{Lebesgue measure} on $\mathbb{T}$ . The norm induced by the inner product is given by 
	 \[
	 \|f\|_{L^2(E)}=\int_{\bbT}\|f\|_{E}^2dm.
	 \]
	 If $\dim E=1$ (i.e., $E = \mathbb{C}$) then $L^2 (E)$ consists of scalar-valued functions and is
	 denoted by $L^2$.
	 
	 The Hardy space $H^2(E)$ is the subspace of $L^2(E)$ formed by the functions with vanishing negative Fourier coefficients; it can be identified with a space of $E$- valued functions analytic in $\mathbb{D}$, from which the boundary values can be recovered almost everywhere through radial limits. One can also view $H^2(E)$ as the direct sum of $d$ standard $H^2$ spaces. We have the orthogonal decomposition 
	 \[
	 L^2(E)=H^2(E)\oplus (zH^2(E))^*.
	 \]
	 
	 The spaces $ L^\infty(E)\subset L^2(E) $ is formed by the essentially bounded functions with values in $ E $; then $ H^\infty(E)\subset H^2(E) $ are the functions in $ L^\infty(E) $ with vanishing negative Fourier coefficients.
	 
	 Taking into account that $\mathcal{L}(E)$ is a Hilbert space endowed with the Hilbert-Schmidt norm, we may similarly define $ H^2(\LL(E)) \subset L^2(\LL(E))$ and $ H^\infty(\LL(E)) \subset L^\infty(\LL(E)) $. Note, however, that we prefer to consider on $ L^\infty(\LL(E)) $ and $ H^\infty(\LL(E)) $ the equivalent norm obtained by considering on $ \LL(E) $ the operator norm instead of the Hilbert-Schmidt norm. 
	 
	 The space $L^2(\mathcal{L}(E))$ may be identified with the matrices with all the entries in $L^2$.
	 We have an orthogonal decomposition
	 \[
	 L^2 (\mathcal{L}(E)) =[zH^2 \mathcal{L}(E)]^*\oplus H^2(\mathcal{L}(E)).\]

	 The space
	 $ L^\infty(\LL(E)) $
	 acts on $L^2 (E)$ by means of multiplication: to $\Phi\in\ L^\infty(\mathcal{L}(E))$ we associate the operator $M_\Phi$ defined by $M_\Phi(f)=\Phi f$ for all $f\in L^2(E)$.

	 Let $S$ denote the forward shift operator $(Sf)( z ) = zf (z)$ on $H^2 (E)$ ; it is the
	 restriction of $M_z$ to $H^2(E)$. Its adjoint (the backward shift) is the operator 
	 \[ 
	 (S^*f)(z) =
	 \frac{f(z )-f(0)}{z}
	 \]
	 One sees easily that $I-SS^*$ is precisely the orthogonal projection onto the space of
	 constant functions. 
	 
	 The main object of study is formed by the \emph{model spaces} and the operators acting on them. These are defined as follows.
	 First, an inner function is an element $ \Theta\in H^2(\mathcal{L}(E))$ whose boundary values are
	 almost everywhere unitary operators in $\mathcal{L}(E)$. The inner function is called \emph{pure} if $ \|\Theta(0)\|<1 $; this is equivalent to requiring that $ \Theta $ has no constant unitary part. 
	 
	 Consider then a pure inner
	 function $\Theta$, with values in $\mathcal{L}(E)$. The model
	 space associated to a pure inner function $\Theta$, denoted by
	 $\mathcal{K}_{\Theta}$ and is defined by
	 $$\mathcal{K}_{\Theta}=H^{2}(E)\ominus \Theta
	 H^{2}(E).$$ \par Just like the Beurling-type subspaces $\Theta
	 H^{2}(E)$ constitute nontrivial invariant subspaces for the
	 unilateral shift $S$ on $H^{2}(E)$, the subspaces
	 $\mathcal{K}_{\Theta}$ play an analogous role for the backward
	 shift $S^{*}$. The orthogonal projection onto
	 $\mathcal{K}_{\Theta}$ will be denoted by $P_{\Theta}$. 
	 It is also known that $ \KK_\Theta\cap H^\infty(E) $ is dense in $ \KK_\Theta $.
	 
	 The analogous space of matrix-valued functions is denoted
	 by $\mathcal{M}_\Theta$; it is the orthogonal complement of $\Theta H^2(\mathcal{L}(E))$ in $H^2(\mathcal{L}(E))$.

	 The model operator $S_\Theta\in\mathcal{L}(\mathcal{K}_\Theta)$ is defined by the formula 
	 \[
	 (S_\Theta f)(z)=P_\Theta (zf).
	 \]
	 The adjoint of $S_\Theta$ is given by
	 \[
	 (S_\Theta^*f)(z)=\frac{f(z)-f(0)}{z};
	 \]
	 it is the restriction of the left shift in $H^2(E)$ to the $S^*$
	 -invariant subspace $\mathcal{K}_\Theta$.

	 Let us assume that  $ \Theta(0)=0 $, so $ \Theta=z\Theta_1 $, which is the case we will consider in the sequel.
	 Then $ I-S_\Theta S_\Theta^* $ is the projection $ P_0 $ onto the constant functions, which are contained in $ \KK_\Theta $, while $ I-S_\Theta^* S_\Theta $ is the projection $P_{D_*}$ onto the space $ \DD_*=\{ \Theta_1 x: x\in E  \} $ (which is also contained in $ \KK_\Theta $).

	 The scalar valued model spaces and operators are obtained when $ \dim E=1 $; that is, when $ E=\bbC $. We have then the classical spaces $ H^2\subset L^2 $ and $L^\infty$. The inner function is a scalar inner function $ \theta $, and the model space is $ \mathcal{K}_{\theta}=H^{2}\ominus \theta
	 H^{2} $. In particular, in case $ \theta(z)=z^n $, $ \KK_\theta $ becomes the $ n $-dimensional space of polynomials of degree at most $ n-1 $.
	 
	 \begin{defn}
	 	A conjugation on a complex Hilbert space $\mathcal{H}$ is a
	 	function $C:\mathcal{H}\longrightarrow \mathcal{H}$ that is
	 	\par
	 	(i) conjugate linear: that is $C(\alpha x+\beta
	 	y)=\overline{\alpha}Cx +\overline{\beta}Cy$ for all $x,y\in
	 	\mathcal{H}$ and
	 	~~~~$\alpha,\beta\in \mathbb{C}$,\par (ii)
	 	involutive: $C^{2}=I$,\par (iii) isometric: $\|Cx\|=\|x\|$ for
	 	all $x\in \mathcal{H}$.
	 \end{defn}
	 
	 The following result is an immediate consequence of a theorem in~\cite{KT}.
	 
	 \begin{lem}\label{le:conjugation} Suppose $ \Theta(0)=0 $, so $ \Theta=z\Theta_1 $. Let $\Gamma$ be a conjugation on $E$, and
	 	suppose that $\Theta(e^{it})^{*}=\Gamma \Theta(e^{it})\Gamma$
	 	a.e. on $\mathbb{T}$. Then
	 	the map $C_{\Gamma}$ defined by
	 	\begin{equation}\label{eq2}
	 	C_\Gamma(f)=z\Theta_{1}\Gamma f
	 	\end{equation}   is a
	 	conjugation on $z\mathcal{K}_{\Theta_1}$.
	 \end{lem}
	 \section{Truncated Toeplitz Operators and Matrix valued Truncated Toeplitz Operators}
	 
	 If $ \phi\in L^\infty $, then the compression of the multiplication operator $ M_\phi $ to $ H^2 $ is called a Toeplitz operator and is denoted by $ T_\phi $. That means that $ T_\phi=P_{H^2}M_\phi|H^2 $, where $P_H^2$ is the orthogonal projection of $L^2$ onto $H^2$. More than a decade ago,
	 Sarason has introduced in~\cite{Sa} the so-called \emph{truncated Toeplitz operators}. Remember that $ P_\theta $ is the orthogonal projection onto the model space $ \KK_\theta $. If $ \phi\in L^\infty $, then the truncated Toeplitz operator $ A^\theta_\phi $ is defined to be the compression of $ M_\phi $ to $ \KK_\theta $. That means that $ A^\theta_\phi=P_\theta M_\phi |\KK_\theta $. In particular, we see that with this notation $ S_\theta= A^\theta_z $.
	 
	 Let us now remember that for $ \theta(z)=z^n $ the space $ \KK_\theta $ is formed by the polynomials of degree not greater than $ n-1 $. The monomials $ z^k $, $ k=0,\dots, n-1 $ form an orthonormal basis of $ \KK_{z^n} $. If we write the matrix of an operator with respect to this basis, then one can see that truncated Toeplitz operators correspond precisely to Toeplitz matrices.
	 
	 Passing now beyond the scalar case, let us suppose that $\Theta$ is a pure inner function. The analogue of truncated Toeplitz operators have been defined in~\cite{KT}, where they are called \emph{matrix valued truncated Toeplitz operators}. 
	 
	 Suppose then that $\Phi\in L^2(\mathcal{L}(E))$. Consider the linear map $f\longrightarrow P_\Theta(\Phi f)$, defined on $\mathcal{K}_\Theta\cap H^\infty(E)$. If it is bounded, then it uniquely determines an operator in $\mathcal{L}(\mathcal{K}_\Theta)$, denoted
	 by $A_\Phi^\Theta$, and called a \emph{matrix valued truncated Toeplitz operator} (MTTO). The function $\Phi$ is then called a \emph{symbol} of the operator. We will usually drop the subscript $\Theta$, as we consider a fixed inner function. We denote by $\mathcal{T}(\mathcal{K}_{\Theta})$ the space of all MTTOs on the model space $\mathcal{K}_\Theta$. 
	 
	 Note that if $ \Phi\in L^\infty(\LL(E)) $ (that is, it is bounded), then it follows that $f\longrightarrow P_\Theta(\Phi f)$ defines a bounded linear operator on the whole of $ \KK_\Theta $, and thus $A_\Phi^\Theta\in \TT(\KK_\Theta)$. But we may have bounded MTTOs which have no bounded symbols, which is one of the complications of the theory. However, a result of \cite{KT} tells that any operator in $ \TT(\KK_\Theta)  $ has a symbol in $ \MM_\Theta+(\MM_\Theta)^* $; this is why we will restrict in the sequel to considering operators $A_\Phi^\Theta  $ with $ \Phi\in \MM_\Theta+(\MM_\Theta)^*  $.\\
	 The operator $S_\Theta$ is a simple example of MTTO; it is obtained by taking $\Phi(z)=zI_E$.  This example is rather special because the symbol is scalar valued.\\
	 It is immediate that 
	 \[
	 A_\Phi^*=A_{\Phi^*};
	 \]
	 so 	$\mathcal{T}(\mathcal{K}_{\Theta})$ is a self adjoint linear space. 
	 
	In section~5 we will see that if $\Theta(z)=z^N I_E$, then the\textsl{} MTTOs obtained are actually familiar objects, namely block Toeplitz matrices of dimension $N$, in which the entries are the matrices of dimension $d$.
	  The theory of Block Toeplitz matrices has been an inspiration for research in the domain of matrix valued truncated Toeplitz operators. In particular, it should be mentioned that some of the classes of block Toeplitz matrices which are closed to multiplication are found in \cite{MAK} and \cite{MAKDT}.\\
	   As supposed above, we will consider a fixed inner function $ \Theta $ and different MTTOs acting on $ \KK_\Theta $. The symbols of these operators will be $ \Phi, \Psi,\dots $.
	  Before the ending of this section we  need to quote result from ~\cite[Chapter VI]{NF}.
	  
	  \begin{prop}\label{p1}
	  	Suppose $\Theta$ be an inner function, and $\Phi\in H^\infty(\mathcal{L}(E))$ such that $\Phi\Theta=\Theta\Phi$. Then
	  	\begin{itemize}
	  		\item[(1)] $\Theta H^2(\mathcal{L}(E))$ is invariant with respect to $M_\Phi$ (and consequently $\mathcal{K}_\Theta)$ is invariant under $M_{\Phi}^*$.
	  		\item [(2)] $A_\Phi S_\Theta=S_\Theta A_\Phi$ and consequently $A_\Phi^*S_\Theta^*=S_\Theta^*$
	  	\end{itemize}
	  \end{prop}
	  Set $\Delta=I-S_\Theta S_\Theta^*$. Note that $P_0$ is the orthogonal projection onto the constant functions contained in $\mathcal{K}_\Theta$.\\
	  
	  The next result from~\cite{KT} characterizes MTTOs among all  operators on $ \KK_\Theta $.
	  
	  \begin{prop}\label{pr1}
	  	A bounded operator $A$ on $\mathcal{K}_\Theta$ belongs to $\mathcal{T}(\mathcal{K}_\Theta)$ if and only if  there exist operators $B,B^\prime$ on $\mathcal{K}_\Theta$ such that 
	  	\[
	  	\Delta(A)=BP_0+P_0{B^\prime}
	  	\]
	  	In this case $A=A_{\Phi+{\Phi^\prime}^*}$, where $\Phi,\Phi^\prime\in H^2(\mathcal{L}(E))$.
	  \end{prop}
	 \section{Main results}
	 In view of the result of Sedlock ~\cite{NS}, a natural question is to determine when is the product of two MTTOs still an MTTO. However, there is no such simple result in the general case, and we need some supplementary assumptions to obtain the main result, Theorem~\ref{th:commutationMTTO}. The path we take is suggested by~\cite{CT}, but the matrix valued situation is much more complicated.

	 We will consider in the rest of this section a fixed inner function $\Theta\in H^\infty(E)$  subjected to the condition $\Theta(0)=0$. Then $\Theta(z)=z\Theta_{1}(z)$, where $\Theta_{1}\in H^\infty(E)$ is also inner. We have the orthogonal decomposition
	 \begin{equation}\label{eq:decomposition Theta1}
	 \KK_\Theta=E\oplus z\KK_{\Theta_1}.
	 \end{equation}
	 
	 Take now $\Phi\in \MM_\Theta+(\MM_\Theta)^* $. We can write then 
	 \begin{equation}\label{eq1}
	 \Phi=z\Phi_{+}+\bar{z}\Phi_{-}^*+\Phi_0
	 \end{equation}
	 with $\Phi_\pm\in\mathcal{M}_{\Theta_1}$ and $\Phi_0\in \mathcal{L}(E)$. 
	 If $ \Phi(e^{it})=\dss\sum_{n=-\infty}^{\infty}\Phi_n e^{int}$ with $\Phi_n\in\mathcal{L}(E)$, then 
	 \begin{equation}\label{eq:formula for Phi+}
	 \begin{split}
	 \Phi_+(z)&=\sum_{n=1}^{\infty}\Phi_n z^n=
	 \sum_{n=1}^{\infty} \left( \int \Phi(e^{it})e^{-int}\, dt \right) z^n\\
	 &=\int \Phi(e^{it}) \left( \sum_{n=1}^{\infty} e^{-int}z^n \right)dt =\int \Phi(e^{it}) \frac{e^{it}z}{1-e^{it}z}dt.
	 \end{split}
	 \end{equation}

	 Remember that two operators $ A,B $ are said to \emph{doubly commute} if $ AB=BA $ and $ AB^*=B^*A $ (whence it follows that also $ A^*B^*=B^*A^* $ and $ A^*B=BA^* $).
	 
	 \begin{lem}\label{le:different commutations}
	 	Suppose that $\Theta(0)=0$ and $\Phi,\Psi\in\mathcal{M}_\Theta+(\MM_\Theta)^*$ such that $\Phi(e^{it})\Psi(e^{is})=\Psi(e^{is})\Phi(e^{it})$ for any $ t,s $. 
	 	\begin{itemize}
	 		\item[(i)] For any $ s,t $ we have  $\Phi_+(e^{it})\Psi(e^{is})=\Psi(e^{is})\Phi_+(e^{it})$  and $\Phi_-(e^{it})\Psi(e^{is})=\Psi(e^{is})\Phi_-(e^{it})$.
	 		
	 		\item[(ii)] If the values of $\Phi,\Psi$ doubly commute with those of $\Theta$, then the same is true for $ \Phi_{\pm}, \Psi_{\pm} $.
	 		
	 		\item[(iii)] If $\Gamma$ is a conjugation on $E$ such that $\Phi(e^{it})^*=\Gamma\Phi(e^{it})\Gamma$, then $\Phi_\pm(e^{it})^*=\Gamma\Phi_\pm(e^{it})\Gamma$.
	 	\end{itemize}
	 \end{lem}
	 
	 \begin{proof}
	 	We will give the proof only for one of the equalities in (ii); the rest are similar.
	 	
	 	Using~\eqref{eq:formula for Phi+}, we have	
	 	\[
	 	\begin{split}
	 	\Psi(e^{is}) \Phi_+(z)&= \Psi(e^{is})\int \Phi(e^{it}) \frac{e^{it}z}{1-e^{it}z}dt=
	 	\int \Psi(e^{is}) \Phi(e^{it}) \frac{e^{it}z}{1-e^{it}z}dt\\
	 	&=\int \Phi(e^{it})\Psi(e^{is}) \frac{e^{it}z}{1-e^{it}z}dt=
	 	\left(\int \Phi(e^{it}) \frac{e^{it}z}{1-e^{it}z}dt\right)\Psi(e^{is})=\Phi_+(z)\Psi(e^{is}).
	 	\end{split}
	 	\]	
	 	By taking radial limits a.e., one obtains the required commutativity.
	 \end{proof}

	 The next lemma gives an identification of elements in $ \MM_\Theta $.
	 
	 \begin{lem}\label{le:MM_Theta}
	 	The map $ \Phi\mapsto J_\Phi $, defined by 
	 	\begin{equation}\label{eq4}
	 	J_\Phi(x)(z)=\Phi(z)x.
	 	\end{equation}
	 	is a bijection between $ \MM_\Theta $ and $ \LL(E, \KK_\Theta) $.
	 \end{lem}
	 
	 \begin{proof}
	 	Fixing a basis $e_1,\cdots ,e_d$ in $E$ and defining the transformation $J:E\longrightarrow\mathcal{K}_\Theta$ as follow
	 	\begin{equation}\label{eq3}
	 	J(e_k)=\phi_k,\qquad 1\leq k\leq d,
	 	\end{equation}
	 	where $\phi_k\in\mathcal{K}_\Theta$ for every $1\leq k\leq d$. If we arrange $\phi_k$ as a column vectors then we obtain $\Phi\in\mathcal{M}_\Theta$.  Conversely if we have $\Phi\in\mathcal{M}_\Theta$ then we obtain the map $J_\Phi: E\longrightarrow\mathcal{K}_\Theta$ which sends $e_k$ to $k$th column of $\Phi$ for every $1\leq k\leq d$. 
	 	
	 	Note that, if $\Phi\in\mathcal{M}_\Theta$, then the relation between $J_\Phi$ and $\Phi$ is simply 
	 	\begin{equation*}\label{eq4}
	 	J_\Phi(x)(z)=\Phi(z)x.\qedhere
	 	\end{equation*}	
	 \end{proof}
	 We denote by $J_0$ the embedding of $E$ in $K_\Theta$; that is for every $x\in E$, $J_0(x)=x$. 
	 It is easy to see that in a given basis the matrices of
	 functions in $z\mathcal{M}_{\Theta_1}$ are characterized by the fact that  columns are functions in $z\mathcal{K}_{\Theta_1}$.
	 
	 Finally, we define $\mathbf{C}_\Gamma(\Phi) $ by giving the action of $J_\Phi$ on $x\in E$ as 
	 \begin{equation}\label{eq:bold C Gamma}
	 J_{\mathbf{C}_\Gamma(\Phi)} x=C_\Gamma(\Phi\Gamma x),
	 \end{equation}
	 where $ C_\Gamma $ is defined by~\eqref{eq2}.

	 In the rest of this section we assume that $\mathfrak{F}$ is a commutative algebra of functions contained in $\MM_\Theta+(\MM_\Theta)^*$, such that all the elements of $\mathfrak{F}$ doubly commute with those of $\Theta$. The next lemma is the main technical result of this section.
	 
	 \begin{lem}\label{l1}
	 	
	 	Suppose that $\Theta(0)=0$ and $\Phi,\Psi\in\mathfrak{F}$. Then there exist operators $ X, Y\in\mathcal{M}_\Theta $ such that 
	 	\[
	 	\Delta(A_\Phi A_\Psi)=J_{z\Phi_+}J_{z\Psi_{-}}^*-J_{\mathbf{C}_\Gamma(z\Phi_{-})}J_{\mathbf{C}_\Gamma(z\Psi_{+})}^*+XP_0+P_0Y
	 	\]
	 \end{lem}

	 \begin{proof}
	 	
	 	For any $\Phi\in\mathcal{M}_{\Theta}+(\MM_\Theta)^*$ we will denote $\hat{\Phi}=\Phi-\Phi_0$.
	 	It follows easily from Lemma~\ref{le:different commutations} that $\Phi_0\Psi(e^{is})=\Psi(e^{is})\Phi_0$
	 	and $\Phi_0\hat{\Psi}(e^{is})=\hat{\Psi}(e^{is})\Phi_0$. In the same way one can obtain $\Psi_0\Phi(e^{is})=\Phi(e^{is})\Psi_0$ and $\Psi_0\hat{\Phi}(e^{is})=\hat{\Phi}(e^{is})\Psi_0$. A similar argument works for double commutation with $\Theta$.\\
	 	Since $\Phi_0,\Psi_0\in H^2(\mathcal{L}(E))$ commutes with $\Theta$ then by using Proposition \ref{p1} $S_\Theta$ commutes with $A_{\Phi_0}$ and $A_{\Psi_0}$, and therefore 
	 	\begin{align*}
	 	\Delta(A_{\Phi} A_{\Psi})
	 	&=\Delta(A_{\hat{\Phi}}A_{\hat{\Psi}})+\Delta(A_{\hat{\Phi}}A_{\Psi_0})+\Delta(A_{\Phi_0}A_{\hat{\Psi}})+\Delta(A_{\Phi_0}A_{\Psi_0}).\\
	 	&=\Delta(A_{\hat{\Phi}}A_{\hat{\Psi}})+\Psi_0\Delta(A_{\hat{\Phi}})+\Phi_0\Delta(A_{\hat{\Psi}})+\Phi_0\Psi_0\Delta(I)\\
	 	&=\Delta(A_{\hat{\Phi}}A_{\hat{\Psi}}) +\Psi_0\Delta(A_{\hat{\Phi}})+\Phi_0\Delta(A_{\hat{\Psi}})+\Phi_0\Psi_0P_0.
	 	\end{align*}
	 	
	 	By Proposition \ref{pr1}  there exist operators $B=z\Phi_{+}$ and $B^\prime=z\Phi_{-}$ such that $\Delta(A_{\hat{\Phi}})=z\Phi_{+}P_0+P_0\bar{z}\Phi_{-}^*$ with $\Phi_{\pm}\in\mathcal{M}_{\Theta_{1}}$.
	 	Similarly $\Delta(A_{\hat{\Psi}})=z\Psi_{+}P_0+P_0\bar{z}\Psi_{-}^*$ and $\Psi_{\pm}\in\mathcal{M}_{\Theta_{1}}$.
	 	
	 	Using Lemma~\ref{le:different commutations}, we have 
	 	\begin{align*}
	 	\Delta(A_\Phi A_\Psi)
	 	&=\Delta(A_{\hat{\Phi}}A_{\hat{\Psi}})+\Psi_0(z\Phi_{+}P_0+P_0\bar{z}\Phi_{-}^*)+\Phi_0(z\Psi_{+}P_0+P_0\bar{z}\Psi_{-}^*)+\Phi_0\Psi_0P_0\\
	 	&=\Delta(A_{\hat{\Phi}}A_{\hat{\Psi}})+(\Psi_0 z\Phi_{+}+\Phi_0z\Psi_{+}+\Phi_0 \Psi_0)P_0+(\Psi_0P_0 \bar{z}\Phi_{-}^*+\Phi_0P_0\bar{z}\Psi_{-}^*)
	 	\end{align*}
	 	Since $P_0$ is the projection onto the constants then it must commute with $\Phi_0$ and $\Psi_0$. Therefore
	 	\begin{equation}\label{eq6}
	 	\Delta(A_\Phi A_\Psi)=\Delta(A_{\hat{\Phi}}A_{\hat{\Psi}})+(\Psi_0 z\Phi_{+}+\Phi_0z\Psi_{+}+\Phi_0\Psi_0)P_0+P_0(\Psi_0\bar{z}\Phi_{-}^*+\Phi_0\bar{z}\Psi_{-}^*)
	 	\end{equation}
	 	Now, by using the definition of $\Delta$,
	 	\begin{align*}
	 	\Delta(A_{\hat{\Phi}}A_{\hat{\Psi}})
	 	&=A_{\hat{\Phi}}A_{\hat{\Psi}}-S_\Theta A_{\hat{\Phi}}A_{\hat{\Psi}}S_\Theta^*\\
	 	&=A_{\hat{\Phi}}A_{\hat{\Psi}}-A_{\hat{\Phi}}S_\Theta A_{\hat{\Psi}}S_\Theta^*+A_{\hat{\Phi}}S_\Theta A_{\hat{\Psi}}S_\Theta^*-S_\Theta A_{\hat{\Phi}}A_{\hat{\Psi}}S_\Theta^*\\
	 	&=A_{\hat{\Phi}}\Delta(A_{\hat{\Psi}})+\Delta(A_{\hat{\Phi}})S_\Theta A_{\hat{\Psi}}S_\Theta^*-S_\Theta A_{\hat{\Phi}}P_{{\mathcal{D}}_*}A_{\hat{\Psi}}S_\Theta^*\\
	 	&=A_{\hat{\Psi}}(z\Phi_{+}P_0+P_0\bar{z}\Psi_{-}^*)+(z\Phi_{+}P_0+P_0\bar{z}\Phi_{-}^*)S_\Theta A_{\hat{\Psi}}S_{\Theta}^{*}-S_\Theta A_{\hat{\Phi}}P_{{\mathcal{D}}_*}A_{\hat{\Psi}}S_\Theta^*
	 	\end{align*}
	 	or
	 	\begin{equation}\label{eq7}
	 	\Delta(A_{\hat{\Phi}}A_{\hat{\Psi}})
	 	=	A_{\hat{\Phi}}z\Psi_{+}P_0+A_{\hat{\Phi}}P_0\bar{z}\Psi_{-}^*+z\Phi_{+}P_0 S_\Theta A_{\hat{\Psi}}S_{\Theta}^{*} +P_0\bar{z}\Phi_{-}^*S_\Theta A_{\hat{\Psi}}S_{\Theta}^{*}-S_\Theta A_{\hat{\Phi}}P_{{\mathcal{D}}_*}A_{\hat{\Psi}}S_\Theta^*.
	 	\end{equation}
	 	Since the constant functions are in $\mathcal{K}_{\Theta}$ , we have $ P_0P_\Theta=P_\Theta P_0 $. Also since $ P_0zf=0 $, then
	 	\[
	 	P_0S_\Theta f=P_0P_\Theta zf= P_\Theta P_0 zf=0,
	 	\]
	 	So the third term in the left hand side of~\eqref{eq7} is 0. 
	 	The second term is 
	 	\[
	 	A_{\hat{\Phi}}P_0\bar{z}\Psi_{-}^*=P_\Theta (z\Phi_+ +\bar z\Phi_-^*)P_0\bar{z}\Psi_{-}^*.
	 	\]
	 	But, since $ \Phi_+\in\mathcal{M}_{\Theta_{1}} $, we have, for any constant function $ x $, $ \Phi_+ x\in\mathcal{K}_{\Theta_1} $ ,   $ z\Phi_+ x\in\mathcal{K}_{\Theta} $, and therefore $ P_\Theta z\Phi_+P_0=z\Phi_+P_0 $. Also,  for any constant function $ x $, $ \bar z\Phi_-^*x\perp H^2(E) $, so $ P_\Theta \bar z\Phi_-^*P_0=0 $. So
	 	\begin{equation}\label{eq8}
	 	\Delta(A_{\hat{\Phi}}A_{\hat{\Psi}})
	 	=	A_{\hat{\Phi}}z\Psi_{+}P_0+z\Phi_+ P_0\bar{z}\Psi_{-}^*+P_0\bar{z}\Phi_{-}^*S_\Theta A_{\hat{\Psi}}S_{\Theta}^{*}-S_\Theta A_{\hat{\Phi}}P_{{\mathcal{D}}_*}A_{\hat{\Psi}}S_\Theta^*.
	 	\end{equation}
	 	Since $ J_0 $ is the embedding of the constants into $ \mathcal{K}_\Theta $, we have, for $ f\in \mathcal{K}_\Theta $ and $ x\in E $,
	 	\[\langle J_0^*f, x\rangle
	 	=\langle f,J_0 x\rangle=
	 	\langle f(0), x\rangle
	 	\] whence $ J^*_0f=f(0) $.
	 	So
	 	\[
	 	J_\Phi J_0^*f=J_\Phi f(0)= \Phi(z)f(0)=\Phi P_0 f,\quad\hbox{for any}\quad f\in \mathcal{K}_\Theta.
	 	\] By taking adjoints we have $J_0 J_\Phi^*=P_0\Phi^*$. Therefore we can write \eqref{eq8} as 
	 	\begin{equation}\label{eq9}
	 	\Delta(A_{\hat{\Phi}}A_{\hat{\Psi}})
	 	=	A_{\hat{\Phi}} J_{z\Psi_{+}}J_0^*+J_{z\Phi_{+}}P_0 J_{z\Psi_{-}}^*+P_0\bar{z}\Phi_{-}^*S_\Theta A_{\hat{\Psi}}S_{\Theta}^{*}-S_\Theta A_{\hat{\Phi}}P_{{\mathcal{D}}_*}A_{\hat{\Psi}}S_\Theta^*.
	 	\end{equation}
	 	Since $ \mathcal{D}_* $ is the space spanned by $ \Theta_1 E $, we can define an isometry $ V:E\to \mathcal{K}_\Theta $ by the formula $ Vx=\Theta_1 x $, and, moreover,  $P_{\mathcal{D}_*}=VV^*  $.
	 	Also, we have
	 	\begin{equation}\label{eq10}
	 	A_{\hat{\Phi}}Vx
	 	=P_\Theta\hat{\Phi}Vx
	 	=P_\Theta z\Phi_+ Vx +P_\Theta \bar z\Phi_-^* Vx.
	 	\end{equation}
	 	Then, using the commutativity between $ \Theta $ and $ \Phi_+ $,
	 	\[
	 	z\Phi_+ V x=
	 	z\Phi_+ \Theta_1 x= z\Theta_1\Phi_+ x=\Theta\Phi_+ x\perp \mathcal{K}_\Theta,
	 	\]
	 	and so the first term in~\eqref{eq10} is 0.\\
	 	We have also 
	 	\begin{align*}
	 	\bar z\Phi_-^*V  x
	 	&=
	 	\bar z\Phi_-^*\Theta_1  x=
	 	\Theta_1 \bar z \Phi_-^* x
	 	= \Theta_1 \bar z \Phi_-^* \Gamma\Gamma x\\
	 	&=\bar z( z\Theta_1 \Gamma (z \Phi_-\Gamma x))
	 	=\bar z C_\Gamma( z \Phi_-\Gamma x).
	 	=\bar{z}J_{\mathbf{C}_\Gamma}(z\Phi_{-})x.
	 	\end{align*}
	 	Since $ C_\Gamma $ is a conjugation on $ z\mathcal{K}_{\Theta_1} $, $ \bar z C_\Gamma (z \Phi_-\Gamma x)\in \mathcal{K}_{\Theta}$, and therefore
	 	\[
	 	P_\Theta\bar z\Phi_-^*V  x=\bar z\Phi_-^*V  x= \bar z C_\Gamma z \Phi_-\Gamma x=\bar{z}J_{\mathbf{C}_\Gamma(z\Phi_{-})}x.
	 	\]
	 	So 
	 	\[
	 	S_\Theta A_{\hat{\Phi}}Vx= z \bar z C_\Gamma z \Phi_-\Gamma x=  C_\Gamma z \Phi_-\Gamma x=J_{\mathbf{C}_\Gamma(z\Phi_{-})}x.
	 	\]
	 	Similarly, we obtain
	 	\[
	 	S_\Theta A_{\hat{\Psi}^*}Vx=C_\Gamma z \Psi_+\Gamma x=J_{\mathbf{C}_\Gamma(z\Psi_{+})} x.
	 	\]
	 	Consequently,
	 	\[
	 	S_\Theta A_{\hat{\Phi}}V=J_{\mathbf{C}_\Gamma(z\Phi_{-})}, \qquad
	 	S_\Theta A_{\hat{\Psi}^*}V=J_{\mathbf{C}_\Gamma(z\Psi_{+})},
	 	\]
	 	Finally, the last term in~\eqref{eq8} is
	 	\begin{equation}\label{eq11}
	 	S_\Theta A_{\hat{\Phi}}P_{{\mathcal{D}}_*}A_{\hat{\Psi}}S_\Theta^*=
	 	S_\Theta A_{\hat{\Phi}}VV^*A_{\hat{\Psi}}S_\Theta^*=
	 	J_{\mathbf{C}_\Gamma(z\Phi_{-})}J_{\mathbf{C}_\Gamma(z\Psi_{+})}^*.
	 	\end{equation}
	 	Combining \eqref{eq9} and \eqref{eq11} we get 	
	 	\[
	 	\Delta(A_{\hat{\Phi}} A_{\hat{\Psi}})=A_{\hat{\Phi}}J_{z\Psi_+}J_0^*+J_{z\Phi_+}P_0J_{z\Psi_{-}}^*+J_0J_{z\Phi_{-}}^*S_\Theta A_{\hat{\Psi}}S_\Theta^*-J_{\mathbf{C}_\Gamma(z\Phi_{-})}J_{\mathbf{C}_\Gamma(z\Psi_{+})}^*
	 	\] 
	 	so we have 
	 	\begin{equation}\label{f_1}
	 	\begin{split}
	 	\Delta(A_\Phi A_\Psi)&=A_{\hat{\Phi}}J_{z\Psi_+}J_0^*+J_{z\Phi_+}P_0J_{z\Psi_{-}}^*+J_0J_{z\Phi_{-}}^*S_\Theta A_{\hat{\Psi}}S_\Theta^*-J_{\mathbf{C}_\Gamma(z\Phi_{-})}J_{\mathbf{C}_\Gamma(z\Psi_{+})}^*\\&+(\Psi_0 z\Phi_{+}+\Phi_0z\Psi_{+}+\Phi_0\Psi_0)P_0+P_0(\Psi_0\bar{z}\Phi_{-}^*+\Phi_0\bar{z}\Psi_{-}^*)
	 	\end{split}
	 	\end{equation}
	 	Now we have the relations
	 	\[
	 	J_\phi J_0^*=\Phi P_0, J_0J_\Phi=P_0\Phi^*
	 	\]
	 	So the second term in \eqref{f_1} becomes 
	 	\[
	 	z\Phi_{+}P_0(z\Phi)^*=J_{z\Phi_{+}}J_0^*J_0J_{z\Phi_{-}}^*=J_{z\Phi_{+}}J_{z\Phi_{-}}^*
	 	\]
	 	Also, since $J_0: E\longrightarrow\mathcal{K}_{\Theta}$ is the embedding of the constants, while $P_0:\mathcal{K}_\Theta\longrightarrow\mathcal{K}_\Theta$ is the projection onto the constant functions, it follows immediately that 
	 	\[
	 	P_0J_0=J_0, J_0^*P_0=J_0^*
	 	\]
	 	Therefore we can write the first and third term in \eqref{f_1} as 
	 	\[
	 	A_{\hat{\Phi}}J_{z\Psi_+}J_0^*=A_{\hat{\Phi}}J_{z\Psi_+}J_0^*P_0,
	 	\quad
	 	J_0J_{z\Phi_{-}}^*S_\Theta A_{\hat{\Psi}}S_\Theta^*=
	 	P_0J_0J_{z\Phi_{-}}^*S_\Theta A_{\hat{\Psi}}S_\Theta^*.
	 	\]
	 	So we have 
	 	\begin{equation}\label{f1}
	 	\begin{split}
	 	\Delta(A_\Phi A_\Psi)&=A_{\hat{\Phi}}J_{z\Psi_+}J_0^*P_0+J_{z\Phi_+}J_{z\Psi_{-}}^*+P_0J_0J_{z\Phi_{-}}^*S_\Theta A_{\hat{\Psi}}S_\Theta^*-J_{\mathbf{C}_\Gamma(z\Phi_{-})}J_{\mathbf{C}_\Gamma(z\Psi_{+})}^*\\&+(\Psi_0 z\Phi_{+}+\Phi_0z\Psi_{+}+\Phi_0\Psi_0)P_0+P_0(\Psi_0\bar{z}\Phi_{-}^*+\Phi_0\bar{z}\Psi_{-}^*)\\
	 	&=J_{z\Phi_+}J_{z\Psi_{-}}^*-J_{\mathbf{C}_\Gamma(z\Phi_{-})}J_{\mathbf{C}_\Gamma(z\Psi_{+})}^*+[A_{\hat{\Phi}}J_{z\Psi_+}J_0^*+\Psi_0 z\Phi_{+}+\Phi_0z\Psi_{+}+\Phi_0\Psi_0]P_0\\&+P_0[J_0J_{z\Phi_{-}}^*S_\Theta A_{\hat{\Psi}}S_\Theta^*+\Psi_0\bar{z}\Phi_{-}^*+\Phi_0\bar{z}\Psi_{-}^*]\\
	 	&=J_{z\Phi_+}J_{z\Psi_{-}}^*-J_{\mathbf{C}_\Gamma(z\Phi_{-})}J_{\mathbf{C}_\Gamma(z\Psi_{+})}^*+XP_0+P_0Y
	 	\end{split}
	 	\end{equation}
	 	where,  $X=A_{\hat{\Phi}}J_{z\Psi_+}J_0^*+\Psi_0 z\Phi_{+}+\Phi_0z\Psi_{+}+\Phi_0\Psi_0$ and $Y=J_0J_{z\Phi_{-}}^*S_\Theta A_{\hat{\Psi}}S_\Theta^*+\Psi_0\bar{z}\Phi_{-}^*+\Phi_0\bar{z}\Psi_{-}^*$.
	 \end{proof}
	 \begin{thm}\label{th:commutationMTTO}
	 	Suppose $\Theta(0)=0$ and $\Phi,\Psi,\chi , \zeta \in\mathfrak{F}$. Then $A_\Phi A_\Psi-A_\chi A_\zeta\in\mathcal{T}(\mathcal{K}_\Theta)$ if and only if
	 	\[	J_{z\Phi_{+}}J_{z\Psi_{-}}^*-J_{\mathbf{C}_\Gamma(z\Phi_-)}J_{\mathbf{C}_\Gamma(z\Psi_{+})}^*=J_{z\chi_+}J_{z\zeta_-}^*-J_{\mathbf{C}_\Gamma(z\chi_-)}J_{\mathbf{C}_\Gamma(z\zeta_+)}^*.
	 	\]
	 \end{thm}
	 \begin{proof}
	 	By Lemma \ref{l1} there exists operators $X, Y\in\mathcal{M}_\Theta$ such that 
	 	\[
	 	\Delta(A_\Phi A_\Psi-A_\chi A_\zeta)=J_{z\Phi_{+}}J_{z\Psi_{-}}^*-J_{\mathbf{C}_\Gamma (z\Phi_{-})}J_{\mathbf{C}_\Gamma(z\Psi_{+})}-J_{z\chi_+}J_{z\zeta_-}^*+J_{\mathbf{C}_\Gamma(z\chi_-)}J_{\mathbf{C}_\Gamma(z\zeta_+)}+XP_0+P_0Y
	 	.\]
	 	By Proposition \ref{pr1}, we have $A_\Phi A_\Psi-A_\chi A_\zeta\in\mathcal{T}(\mathcal{K}_\Theta)$ if and only  there exist operators $B, B^\prime$ such that
	 	\begin{equation*}\label{eq:condition}
	 	\Delta(A_\Phi A_\Psi-A_\chi A_\zeta)=BP_0+P_0B.
	 	\end{equation*}
	 	The last two equations say that $A_\Phi A_\Psi-A_\chi A_\zeta\in\mathcal{T}(\mathcal{K}_\Theta)$ if and only if there exist $ X',Y' $ such that
	 	\begin{equation}\label{sum}
	 	J_{z\Phi_{+}}J_{z\Psi_{-}}^*-J_{\mathbf{C}_\Gamma (z\Phi_{-})}J_{\mathbf{C}_\Gamma(z\Psi_{+})}-J_{z\chi_+}J_{z\zeta_-}^*+J_{\mathbf{C}_\Gamma(z\chi_-)}J_{\mathbf{C}_\Gamma(z\zeta_+)}=X'P_0+P_0Y^\prime.
	 	\end{equation}
	 	Since $\Theta(0)=0$ then we can write  $\mathcal{K}_\Theta=E\oplus\mathcal{K}_{\Theta_{1}}$. With respect to this decomposition, the left hand side of \eqref{sum} has zeros on the first row and column, while the right hand side is the general form of an operator that has zeros in the lower right corner. Now it is clear that for~\eqref{sum} to be true both sides have to be zero, which proves the theorem. 
	 \end{proof}
 
	 The following result is the main result of this paper: it gives the answer to the question stated at the beginning of this section, namelu when is the product of two MTTOs also an MTTO. 
	 \begin{thm}\label{th:multiplicationMTTO}
	 	Suppose $\Theta(0)=0$, $\Phi,\Psi\in\mathfrak{F}$ and $A_\Phi, A_\Psi, \in\mathcal{T}(\mathcal{K}_\Theta)$. Then $A_\Phi A_\Psi\in\mathcal{T}(\mathcal{K}_\Theta)$ if and only if $J_{z\Phi_{+}}J_{z\Psi_{-}}^*=J_{\mathbf{C}_\Gamma(z\Phi_-)}J_{\mathbf{C}_\Gamma(z\Psi_{+})}^*.$
	 \end{thm}
 \begin{proof}
 	Applying Theorem \ref{th:multiplicationMTTO} to the case $\chi=\zeta=0$.
 \end{proof}
 \section{A particular case: Block Toeplitz Matrices}
 In this section let $\Theta(z)=z^NI_E$ for some fixed positive integer $N$. Then $\KK_\Theta$ is the Hilbert space of all polynomials in $z$ of degree at most $N-1$ with coefficients from $E$,i.e.,
 \[
 \KK_\Theta=\{a_0+a_1z+\cdots a_{N-1}z^{N-1};a_0,a_1,\cdots,a_{N-1}\in E\}.
 \]
 One can also identify this space with N-copies of $E$ by mapping $\dss\sum_{k=0}^{N-1}a_kz^k$ into $\dss\otimes_{k=0}^{N-1}a_k$. Now if we take $A_\Phi$ corresponding to function $\Phi\in L^{2}(\LL(E))$ having Fourier expansion $\Phi(e^{it})=\dss\sum_{n=-\infty}^{\infty}\Phi_{n}e^{int}$, with $\Phi_n\in\LL(E)$. Then it can be easily seen that with respect to the direct decomposition given above $A_\Phi$ has a natural representation as a matrices: Let  $\Phi({k})$ be the $k$-th Fourier coefficient of $\Phi$. Then 
 $$\TT(\KK_\Theta)=\Biggl\{A_\Phi=\begin{pmatrix}
 	\Phi(0) & \Phi(1)&\cdots &\Phi(N-1)\\
 	\Phi(-1)& \Phi(0)&\cdots& \Phi(N-2)\\
 	\vdots & \vdots & \ddots&\vdots \\
 	\Phi(1-N)& \Phi(2-N)& \cdots &\Phi(0)
 \end{pmatrix}; \Phi({k})\in\LL(E)
 \Biggr\}
 $$
 Suppose now that $\Phi,\Psi$ belong to a commutative algebra $\mathfrak{F}$. Since $\Theta$ is a scalar valued inner function the double commutation assumption on $\Phi,\Psi$ is satisfied. As stated in section~4, we can write the symbol $\Phi\in\mathfrak{F}$ of any MTTO as 
 \[
 \Phi=z\Phi_{+}+\bar{z}\Phi_{-}+\Phi_{0},
 \]
 where $\Phi_{\pm}\in\MM_{\Theta_1}$ and $\Phi_{0}\in\LL(E)$. Let $\Phi_+(k)$ and $\Phi_{-}(k)$ denote the Fourier coefficients of $\Phi_{+}$ and $\Phi_{-}$ respectively then  we have  
 \[
 A_{z\Phi_{+}}=\begin{pmatrix}
 0 & 0 &\cdots & 0\\
 \Phi_+(0) & 0&\cdots& 0  \\
 \vdots & \vdots & \ddots&\vdots \\
 \Phi_+(N-2) & \Phi_+(N-3)& \cdots &0
 \end{pmatrix}
 \] 
 and 
  \[
 A_{\bar{z}\Phi_{-}^{*}}=\begin{pmatrix}
 0 & \Phi_{-}^{*}(0) &\cdots & \Phi_{-}^{*}(N-2) \\
 0 & 0&\cdots& \Phi_{-}^{*}(N-3)   \\
 \vdots & \vdots & \ddots&\vdots \\
 0 & 0 & \cdots &0
 \end{pmatrix}
 \]
 We have therefore 
 \[
 A_\Phi=A_{z\Phi_{+}+\bar{z}\Phi_{-}+\Phi_{0}}=
 \begin{pmatrix}
 \Phi_{0} & \Phi_{-}^{*}(0) &\cdots & \Phi_{-}^{*}(N-2) \\
 \Phi_{+}(0) & \Phi_{0}&\cdots& \Phi_{-}^{*}(N-3)   \\
 \vdots & \vdots & \ddots&\vdots \\
 \Phi_{+}(N-2) & \Phi_{+}(N-3) & \cdots &\Phi_{0}
 \end{pmatrix}
 \] 
 Similarly, 
  \[
 A_\Psi=A_{z\Psi_{+}+\bar{z}\Psi_{-}+\Psi_{0}}=
 \begin{pmatrix}
 \Psi_{0} & \Psi_{-}^{*}(0) &\cdots & \Psi_{-}^{*}(N-2) \\
 \Psi_{+}(0) & \Psi_{0}&\cdots& \Psi_{-}^{*}(N-3)   \\
 \vdots & \vdots & \ddots&\vdots \\
 \Psi_{+}(N-2) & \Psi_{+}(N-3) & \cdots &\Psi_{0}
 \end{pmatrix}
 \] 
Since $\Theta(z)=z^NI_E$, 
$\MM_\Theta$ is the space of all polynomials in $z$ of degree at most $N-1$ with coefficients from $\LL(E)$; also, $\Theta(0)=0$ imply that $\Theta(z)=z\Theta_1(z)$, where $\Theta_1$ is also inner. 
Then $\MM_{\Theta_1}$ is also the space of polynomials in $z$ with coefficients from $\LL(E)$ but of degree not greater than $N-2$.
For any $x\in E$, $J_{z\Phi_{+}}x=z\Phi_{+}x$, where $\Phi_{+}\in\LL(E)$, and for any $f\in\KK_\Theta$, $J_{z\Psi_{-}}^*f=\bar{z}\Psi_{-}^*f$, with $\Psi_{-}\in\LL(E)$.
Now we have 
\begin{align*}
 J_{z\Phi_{+}}J_{z\Psi_{-}}^*f
 &=z\Phi_{+}\bar{z}\Psi_{-}^*f=\Phi_{+}\Psi_{-}^*f=\left(\sum_{k=0}^{N-2}\Phi_{+}(k)z^k\right)\left(\sum_{k=0}^{N-2}\Psi_{-}^*(k)\bar{z}^k\right)f\\
 &=\left(\sum_{k=0}^{N-2}\Phi_{+}(k)\Psi_{-}^*(k)+\sum_{k=0}^{N-3}\Phi_{+}(k)\Psi_{-}^*(k+1)\bar{z}+\cdots\Phi_{+}(0)\Psi_{-}^*(N-2)\bar{z}^{N-2}\right)f+\\
 &\left(\sum_{k=0}^{N-3}\Phi_{+}(k+1)\Psi_{-}^*(k)z+\sum_{k=0}^{N-4}\Phi_{+}(k+2)\Psi_{-}^*(k)z^2+\cdots\Phi_{+}(N-2)\Psi_{-}^*(0)z^{N-2}\right)f\\
 &=\sum_{m=N-2}^{0}\sum_{k=0}^{m}\Phi_{+}(k)\Psi_{-}^*(k+(N-2)-m)\bar{z}^{N-2-m}f+\\
 &\sum_{m=N-3}^{0}\sum_{k=0}^{m}\Phi_{+}(k-m+N-2)\Psi_{-}^*(k)z^{N-2-m}f
 \end{align*} 
\begin{equation}\label{multilpcation1}
\begin{split}
J_{z\Phi_{+}}J_{z\Psi_{-}}^*f=
\sum_{m=N-2}^{0}\sum_{k=0}^{m}\Phi_{+}(k)\Psi_{-}^*(k+(N-2)-m)\bar{z}^{N-2-m}f+\\
+\sum_{m=N-3}^{0}\sum_{k=0}^{m}\Phi_{+}(k+(N-2)-m)\Psi_{-}^*(k)z^{N-2-m}f
\end{split}
\end{equation}
and for any $f$ in $\KK_\Theta$ \begin{equation}\label{multilpcation2}
\begin{split}
J_{\mathbf{C}_\Gamma(z\Phi_{-})} J_{\mathbf{C}_\Gamma(z\Psi_{+})}^*f=\Phi_{-}^*\Psi_{+}f
=\sum_{m=N-2}^{0}\sum_{k=0}^{m}\Phi_{-}^*(k)\Psi_{+}(k+(N-2)-m)\bar{z}^{N-2-m}f+\\
+\sum_{m=N-3}^{0}\sum_{k=0}^{m}\Phi_{-}^*(k-m+N-2)\Psi_{+}(k)z^{N-2-m}f
\end{split}
\end{equation}
By Theorem~ \ref{th:multiplicationMTTO} $A_\Phi A_\Psi\in\TT(\KK_\Theta)$ if and only if  $J_{z\Phi_{+}}J_{z\Psi_{-}}^*=J_{\mathbf{C}_\Gamma(z\Phi_{-})} J_{\mathbf{C}_\Gamma(z\Psi_{+})}^*$. Comparing corresponding coefficients we get   
\begin{equation}\label{eq0}
\bar{z}^{N-2}:\Phi_{+}(0)\Psi_{-}^*(N-2)=\Phi_{-}(0)^*\Psi_{+}(N-2)\implies \Phi_{+}^*(0)\Psi_{-}(N-2)=\Phi_{-}^*(0)\Psi_{+}(N-2)\\
\end{equation}
\begin{equation}\label{eq00}
\bar{z}^{N-3}:\Phi_{+}(0)\Psi_{-}^*(N-2)+\Phi_{+}(1)\Psi_{-}^*(N-3)=\Phi_{-}^*(0)\Psi_{+}(N-2)+\Phi_{-}^*(1)\Psi_{+}(N-3)
\end{equation}
Using \eqref{eq0} in \eqref{eq00} we have $\Phi_{+}(1)\Psi_{-}^*(N-3)=\Phi_{-}^*(1)\Psi_{+}(N-3)$. In general 
\begin{equation}\label{Toeplitzproduct1}
\Phi_{+}(i)\Psi_{-}^*(N-2-i)=\Psi_{+}(N-2-i)\Phi_{-}^*(0)\quad \hbox{for every}\quad i=0,\cdots ,N-2. 
\end{equation}
 In the same way compairing coefficients of $z,z^2,\cdots z^{N-2}$ we obtain 
\begin{equation}\label{Toeplitzproduct2}
\Psi_{+}(i)\Phi_{-}^*(N-2-i)=\Phi_{+}(N-2-i)\Psi_{-}^*(i)\quad \hbox{for every}\quad i=0,1,\cdots ,N-2
\end{equation}
If we take $A_i=\Phi_{-}^*(i-1)$, $A_{i-N}=\Phi_{+}(i-(N-1))$ and $B_i=\Psi_{-}^*(i-1)$, $B_{i-N}=\Psi_{+}(i-(N-1))$ for every $i=1,2,\cdots, N-1$ then Lemma 3.1(i) of \cite{MAKDT} imply that $A_\Phi A_{\Psi}$ is a block Toeplitz matrix. Thus the condition  $J_{z\Phi_{+}}J_{z\Psi_{-}}^*=J_{\mathbf{C}_\Gamma(z\Phi_{-})} J_{\mathbf{C}_\Gamma(z\Psi_{+})}^*$ is equivalent to the condition Lemma 3.1(i) of \cite{MAKDT}.

\section*{Acknowledgements}
The author is highly grateful to Dr. Dan Timotin for his valuable suggestions and comments.

\end{document}